\documentclass[12pt]{scrartcl}

\usepackage[
theoremdefs,
biblatex,
final,
]%
{allergy}

\usepackage{unixode}
\usepackage{helvetic}

\usetikzlibrary{matrix}
\tikzset{commdiag/.style={matrix of math nodes, row sep=3.5em, column sep=3em, text height=1.5ex, text depth=0.25ex,ampersand replacement=\&},
smallcommdiag/.style={commdiag, row sep=2em, column sep=1.5em, font=\small},
exseq/.style={commdiag, column sep=2em},
diagequal/.style={double, double distance=2pt, -},
>=stealth,
ineq/.style={baseline=(current  bounding  box.center)},
amph/.style={very thick},
}

\newcommand*\group[1]{\mathrm{#1}}
\newcommand*\aff{\mathsf{a}}
\NewDocumentCommand\Aff{O{d}}{\mathrm{Aff}\paren{#1}}

\newcommand*\Rel[1]{\overline{#1}}
\newcommand*\Hom{\operatorname{Hom}}

\newcommand*\inv{^{-1}}
\newcommand*\RR{\mathbf{R}}
\newcommand*\dd{\mathrm{d}}

\newcommand*\act{\cdot}
\newcommand*\mul{\,}

\newcommand*\Tan{\mathsf{T}}

\newcommand*\Data{\mathcal{D}}
\newcommand*\Comp{\mathcal{C}}
\newcommand*\Alg{{F}}

\NewDocumentCommand\Cont{O{\infty}}{C^{#1}}

\newcommand*\Diff{\mathrm{Diff}}
\newcommand*\affgrp[2]{\begin{bmatrix} #1 & #2\\ 0 & 1 \end{bmatrix}}
\newcommand*\GL[1][d]{\ensuremath{\group{GL}(#1)}}
\NewDocumentCommand\affrel{O{\aff}mm}{#2 \stackrel{#1}\rightsquigarrow{}#3} 

\NewDocumentCommand\AS{O{d}}{\mathbf{A}_{#1}}

\usepackage{authblk}

\author{Olivier Verdier}
\affil{Department of Computing, Mathematics and Physics\\Bergen University College, Bergen, Norway}
\title{Full affine equivariance and weak natural transformations in numerical analysis - the case of B-Series.}

\begin{document}
\maketitle
\begin{abstract}
Many algorithms in numerical analysis are affine equivariant: they are immune to changes of affine coordinates.
This is because those algorithms are defined using affine invariant constructions.
There is, however, a crucial ingredient missing: most algorithms are in fact defined regardless of the underlying dimension.
As a result, they are also invariant with respect to non-invertible affine transformation from spaces of different dimensions.
We formulate this property precisely: these algorithms fall short of being natural transformations between affine functors.
We give a precise definition of what we call a \emph{weak natural transformation} between functors, and illustrate the point using examples coming from numerical analysis, in particular B-Series.
\end{abstract}

\section{Affine Equivariance}
\label{sec-equi}

We define an \emph{algorithm} as a function $\Alg$ from a \emph{data space $\Data$} to a \emph{computation space $\Comp$}:
\begin{equation}
\Alg \colon \Data \to \Comp
.
\end{equation}
In most of the examples, $\Data$ and $\Comp$ are manifolds.

Recall that for a given dimension $d$, the affine group $\Aff$ is defined as the semi-direct product $\Aff \coloneqq \GL \ltimes \RR^d$.
\itodo*{Reference here?}
The action of $\Aff$ on $\RR^d$ is defined as follows.
An element
\begin{align}
g = \affgrp{A}{b}
,
\qquad
A\in\GL
,
\quad
b \in \RR^d
,
\end{align}
acts on an element
\begin{align}
\begin{bmatrix}
x \\ 1	
\end{bmatrix}
\end{align}
by matrix multiplication.
The action is thus
\begin{align}
g \act x = A\mul x + b
.
\end{align}

\begin{definition}
\label{def:algaffequi}
  Suppose that the group $\Aff$ acts on the spaces $\Data$ and $\Comp$.
An algorithm $\Alg \colon \Data \to \Comp$ is \emph{affine equivariant} if the following diagram commutes, for any $\aff \in \Aff$.

\begin{equation}
    \begin{tikzpicture}[ineq]
\matrix(m) [commdiag,]
{
    d_1 \& d_2 \\
    \Alg(d_1) \& \Alg(d_2)\\
};
\path[->]
 (m-1-1) edge node[auto]{$\aff$}     (m-1-2)
 (m-2-1) edge node[auto]{$\aff$}     (m-2-2)
 (m-1-1) edge node[auto] {$\Alg$} (m-2-1)
 (m-1-2) edge node[auto] {$\Alg$} (m-2-2)
;
\end{tikzpicture}
\end{equation}
\end{definition}

In practice, affine equivariance means invariance with respect to a change of affine coordinates.
It means in particular invariance with respect to
\begin{itemize}
    \item translations (change of origin)
    \item anisotropic scalings (change of units)
    \item rotations
    \item shearing
\end{itemize}

We can rephrase \autoref{def:algaffequi} in order to prepare for \autoref{sec-fullaffineequi}.
We regard the group $\Aff$ as a \emph{category} with one object $\star$~\cite[\S\,4.3]{Aw10}.
We regard $\Comp$ and $\Data$ as objects in the category of smooth manifolds and smooth maps.
The actions of $\Aff$ on $\Data$ and $\Comp$ now define \emph{functors}.
With a slight abuse of notation we note these functors $\Data$ and $\Comp$, so the actual data and computational spaces are $\Data_{\star} \coloneqq \Data(\star)$ and $\Comp_{\star}\coloneqq\Comp(\star)$.
Affine equivariance (\autoref{def:algaffequi}) now expresses that the algorithm $\Alg$ is a \emph{natural transformation} from the functor $\Data$, to the functor $\Comp$.
Indeed, \autoref{def:algaffequi} can be rewritten as
\begin{equation}
\label{eq:nattrans}
 \Alg \circ \Data(\aff) = \Comp(\aff) \circ \Alg
\end{equation}
or, using a commuting diagram,
\begin{equation}
    \begin{tikzpicture}[ineq]
\matrix(m) [commdiag,]
{
    \Data_{\star} \& \Data_{\star} \\
    \Comp_{\star} \& \Comp_{\star}\\
};
\path[->]
 (m-1-1) edge node[auto]{$\Data(\aff)$}     (m-1-2)
 (m-2-1) edge node[auto]{$\Comp(\aff)$}     (m-2-2)
 (m-1-1) edge node[auto] {$\Alg$} (m-2-1)
 (m-1-2) edge node[auto] {$\Alg$} (m-2-2)
;
\end{tikzpicture}
.
\end{equation}

We give many examples of such functors $\Data$ and $\Comp$ in this note.

For convenience, for an affine map $\aff(x) = A x + b$, we introduce the corresponding \emph{tangent map}
\begin{equation}
  \Tan \aff \coloneqq A
  .
\end{equation}
Note that the notion of a tangent map is defined for any nonlinear map: indeed, most of what we present in this section can be generalised to any group action (see \autoref{rem-homogeneous}).

\begin{example}[Quadrature]
 The data domain consists of all intervals and continuous functions on those intervals.
This is the union of the spaces $\Cont[0]([\alpha,\beta])$, where we keep track of the interval $[\alpha,\beta]$, so a piece of data is $d = \paren{\alpha,\beta,f}$, where $f\in\Cont[0]([\alpha,\beta])$.
Formally, the data domain $\Data$ is thus a fibre bundle.

The action of $\Aff[1]$ on a pair $(\alpha, \beta)$ is the diagonal action
\begin{equation}
\aff \act (\alpha,\beta) = (\aff(\alpha),\aff(\beta))
\end{equation}
and the action of $\Aff[1]$ on $\Cont[0]([\alpha,\beta])$ is defined by
\begin{equation}
\aff \act f \coloneqq f \circ \aff\inv
,
\end{equation}
so the total action is a fibre bundle mapping (it preserves the fibres).

The action of $\Aff[1]$ on the computational domain $\Comp = \RR$ is the linear action
\begin{equation}
\aff \act x = \Tan\aff \mul x
.
\end{equation}
Affine equivariance is now the requirement that the quadrature formula should fulfil
\begin{equation}
  \label{eq:quadequi}
\Alg\paren[\big]{\aff \cdot (I,f)} = \aff \cdot \Alg(I,f)
.
\end{equation}
Explicitly, this corresponds to the requirement that $\Alg$ should behave as the exact integral under affine transformations.
Indeed, the exact integral fulfils
\begin{equation}
  \int_{a\alpha + b}^{a\beta +b} f\paren[\big]{(y-\beta)/\alpha} \dd y = a \int_{\alpha}^{\beta} f(x) \dd x
\end{equation}
for any real numbers $\alpha \neq 0$ and $\beta$, which is equivalent to~\eqref{eq:quadequi}.

If we interpret the group actions as functors, then the data functor $\Data$ maps the group object $\star$ to $\Data_{\star}$, the fibre bundle defined above, and an invertible one-dimensional affine map $\aff$ is mapped to $\Data(\aff)$ defined by
\begin{equation}
  \Data(\aff)\paren{\alpha,\beta,\Cont[0]([\alpha,\beta],\RR)} = \paren{\aff \mul\alpha,\aff\mul\beta,f\circ \aff\inv}
.
\end{equation}
Similarly, the computational functor $\Comp$ maps the group object $\star$ to $\Comp_{\star} = \RR$, and an invertible one-dimensional affine map $\aff$ is mapped to $\Comp(\aff)$ defined by
\begin{equation}
  \Comp(\aff) \coloneqq \Tan\aff
  .
\end{equation}
\end{example}

\begin{example}[Numerical integrators]
\label{ex:equinuminteg}
    We consider numerical integration of ordinary differential equations (ODEs).
    The data domain $\Data$ is the set of compactly supported vector fields on an affine space $\AS$ of fixed dimension $d$:
    \begin{equation}\Data = \mathfrak{X}_0 (\AS).\end{equation}
    
    The affine action of the group $\Aff[d]$ on a vector field is defined by
    \begin{equation}
    \aff \act f \coloneqq \Tan\aff \circ f \circ \aff\inv
    .
    \end{equation}
    The computational domain $\Comp$ is the set of diffeomorphisms:
    \begin{equation}
    \Comp = \Diff(\AS)
    .
    \end{equation}
    The action on a diffeomorphism is the adjoint action
    \begin{equation}
    \aff \act \Phi \coloneqq \aff \circ \Phi \circ \aff\inv
    .
    \end{equation}
    The equivariance assumption is thus
    \begin{equation}
    \Alg(\aff \act f) = \aff \act \Alg(f)
    .
    \end{equation}
    Again, that requirement makes sense as the exact solution fulfils that property for any invertible mappings (not only the affine ones)~\cite[\S\,2.4]{aromatic}.
    Enforcing equivariance with respect to all invertible transformation would leave us with the exact solution alone.

    For an initial condition $x_0$, the condition means that
    \begin{equation}
    \label{eq:simpleintegequi}
    \Alg(\Tan \aff f \circ \aff\inv)(\aff x_0) = \aff \Alg(f) (x_0)
    .
    \end{equation}
    The layman description of the invariance of an integrator is that if one moves both the initial condition and the vector field with an affine transformation, then the computed point is also moved by the same affine transformation.

    Let us see what that definition becomes for a concrete example of an integrator, the forward Euler method.
    In that case,
    \begin{equation}
    \Alg(f) = \bracket[\big]{x \mapsto x + f(x)}
    .
    \end{equation}

    Writing $\aff\mul x = A x + b$, we can check that
    \begin{align}
\Alg\paren[\Big]{A f \paren[\big]{A\inv(y-\beta)}}(A x_0 + \beta) 
&= (A x_0 + \beta) + A f(x_0) \\ 
&= A (x_0 + f(x_0)) + \beta\\ 
&= \aff(\Alg(f)(x_0))
                                                                                                                                   ,
    \end{align}
    which was condition~\eqref{eq:simpleintegequi}.
    
Using the functor description, the data domain is now $\Data_{\star} = \mathcal{X}_0(\AS)$, and an invertible affine transformation $\aff$ is mapped to $\Data(\aff)\in \Hom(\Data_{\star},\Data_{\star})$ defined by
\begin{equation}
  \Data(\aff)(f) = \Tan\aff \circ f \circ \aff\inv
  .
\end{equation}
The computational functor $\Comp$ maps the group object to $\Comp_{\star} = \Diff(\AS)$, and an invertible map $\aff$ is mapped to $\Comp(\aff)\in\Hom(\Comp_{\star},\Comp_{\star})$ as
\begin{equation}
  \Comp(\aff)(\Phi) = \aff \circ \Phi \circ \aff\inv
  .
\end{equation}

\end{example}

\begin{remark}
  It turns out that all Runge--Kutta methods are affine equivariant.
  It has therefore been conjectured that Runge--Kutta methods, or more precisely, B-Series methods, were the only integrators enjoying that property.
  A recent result shows that this is not the case~\cite{aromatic}.
  An example of an integrator which is affine equivariant but not a B-Series method is
\begin{equation}
  \Alg(f) \coloneqq \bracket[\Big]{x \mapsto f(x)\paren[\big]{1 + \operatorname{div}(f)(x)}}
  .
\end{equation}
See \autoref{ex:bseries} for a complete characterisation of B-Series.
\end{remark}

\begin{example}[Polynomial interpolation and splines]
\label{ex:interpolation}
Here the domain is $\AS^n$, the data of $n$ points $P_i$ in an affine space $\AS$ of dimension $d$.
The computation is a curve $\Cont(\RR,\AS)$, which interpolates the points $P_i$ in a variety of generalised meanings: exact interpolation, splines of various smoothness, etc.

The actions of the affine group $\Aff[d]$ are particularly simple in this case.
As we shall see in \autoref{ex:interpolation2}, this simplicity reflects the fact that the algorithm is in this case a natural transformation between affine functors.

On the domain, $\Data = \AS^n$, the action is the standard diagonal action
\begin{equation}
\aff \act (P_1,\ldots,P_n) \coloneqq (\aff(P_1),\ldots,\aff(P_n))
.
\end{equation}

The action is essentially the same on $\Comp$:
\begin{equation}
\aff \act \gamma \coloneqq \aff \circ \gamma
.
\end{equation}

The equivariance condition is thus simply that
\begin{equation}
\Alg\paren[\big]{\aff\mul P_1,\ldots,\aff \mul P_n} = \aff \mul \Alg(P_1,\ldots,P_n)
.
\end{equation}

The interpretation is particularly intuitive: first moving the control points and then computing the interpolation curves gives the same result as first computing the interpolation curve and then moving it with the same displacement.

Examples of interpolation algorithms which are affine equivariant are
\begin{itemize}
  \item exact interpolation
  \item Bézier splines~\cite[\S\,2.2]{PrBoPa13}
  \item B-Splines~\cite[\S\,5.7]{PrBoPa13}
\end{itemize}

Again, we give the functorial point of view for completeness.
The data functor $\Data$ maps the group object to $\Data_{\star} = \AS^n$, and an invertible affine map $\aff$ is mapped to $\Data(\aff)\in\Hom(\Data_{\star},\Data_{\star})$ defined by $\Data(\aff)(P_i) = (\aff \mul P_i)$. The computational domain is $\Comp_{\star} = \Cont(\RR,\AS)$, and the functor $\Comp$ maps an invertible affine map $\aff$ to $\Comp(\aff)\in\Hom(\Comp_{\star},\Comp_{\star})$ defined by $\Comp(\aff)(\gamma) = \aff\circ\gamma$.

\end{example}

\begin{example}[Downhill simplex minimization algorithm (Nelder--Mead)]
\label{ex:downhill}
The space is $\AS$, and the data is a function $\varphi \in \Cont[0]\paren{\AS}$ to minimize, as well as a set of $n$ starting points.
The algorithm $\Alg$ then produces a new set of $n$ points.

 Usually there are $n = d+1$ points, which span a simplex (hence the name of the algorithm), but this is too restrictive, as we shall see in \autoref{ex:downhill3}.

There are several variants to that algorithm, but the crucial aspect here is that they are all affine invariant~\cite{GoWe10}.

The action on functions is given by $\aff \act \varphi \coloneqq \varphi \circ \aff\inv$, and the action on points is again the diagonal one: $\aff \act x_i = \aff(x_i)$.
The requirement of equivariance is thus $\Alg\paren[\big]{\aff \act (\varphi, x_i)} = \aff \act \Alg(\varphi, x_i)$.
Of course, the actual Nelder--Mead algorithm consists of $N$ iterations of the function $\Alg$ until convergence, and the iterated function $\Alg^N$ inherits the equivariance property of $\Alg$.

What is the meaning of affine equivariance in that case?
It is the idea that if one transforms the function to minimize with an affine transformation, and if one transforms the initial simplex by the same transformation, the final result will be the same as if one had run the algorithm directly, only transforming the last simplex.

The functorial point of view is $\Data_{\star} = \AS^n \times \Cont[0](\RR^d)$, with corresponding action $\Data(\aff)(\varphi,X) = (\varphi\circ \aff\inv, \aff \mul X)$.
The computational domain is $\Comp_{\star} = \AS^n$ and $\Comp(\aff)(X) = \aff \mul X$.
In both cases, $\aff \mul X$ denotes the diagonal action on an element $X\in\AS^n$.
\end{example}

\begin{remark}
\label{rem-homogeneous}
  Even though the main focus of this section is the affine group, it is legitimate to ask which algorithms are equivariant with respect to another group.
Note that \autoref{def:algaffequi} is unchanged: we only replace the affine group with another Lie group, with suitable actions on the data domain $\Data$ and computational domain $\Comp$.

There are already some answers if we restrict the discussion to \emph{numerical integrators on homogeneous spaces}, where equivariance is described along the lines of \autoref{ex:equinuminteg}.
If the underlying homogeneous space is \emph{reductive}, then all the standard extensions of Runge--Kutta methods on homogeneous spaces (Crouch--Grossman, RKMK, commutator-free) are equivariant with respect to the group at hand~\cite{MkVe15}.

When the homogeneous space is symplectic, there are no general way to construct equivariant, symplectic integrators.
Most of such symplectic homogeneous spaces are coadjoint orbits of a Lie group.
In many cases, the construction is still possible using appropriate \emph{symplectic realisations}, i.e., Poisson maps from a symplectic vector space into the Lie--Poisson space at hand.
One thus obtains symplectic integrators on coadjoint orbits, which are symplectic homogeneous spaces.
Those integrators are automatically equivariant with respect to the Lie group at hand~\cite{collective, collectiveHopf}.

\end{remark}

\begin{remark}
  Affine transformations play also a fundamental role in finite element methods.
  The families of polynomial differential forms discovered by Raviart, Thomas, Nédélec, later put in a common framework by Hiptmair, all have a common point: they are all affine invariant spaces: they are mapped to themselves by invertible affine maps (see~\cite[\S\,1.3]{ArFaWi10} and references therein).
  Remarkably, one can describe \emph{all} such spaces~\cite[Th. 3.6]{ArFaWi06}.
  The techniques used are very similar to those used to describe all the affine equivariant integrators in~\cite{aromatic}.
\end{remark}

\section{Full affine equivariance}
\label{sec-fullaffineequi}

In almost all the examples of \autoref{sec-equi}, the algorithms are in fact defined \emph{in any dimension}.
So we have instead a sequence of algorithms $\Alg_d$ for every dimension $d$, mapping a data domain $\Data_d$ into  a computational domain $\Comp_d$.

For instance, an interpolation algorithm is defined for any dimension $d$, takes $n$ points in $\AS$ as input, and returns a curve in $\AS$.
As a result, the data domain is $\Data_d \coloneqq \AS^n$ and the computational domain is $\Comp_d \coloneqq \Cont(\RR,\AS)$.

The crucial observation is that \emph{these functions $\Alg_d$ must be related}.
What we proceed to do now is to express this precisely.
This will sometimes lead to surprising results (see the characterisation of \autoref{ex:bseries}).

We first motivate on an examples why full affine equivariance is needed.
\begin{example}[Downhill simplex minimization algorithm revisited]
\label{ex:downhill2}
We revisit \autoref{ex:downhill}.
For each natural number $d$, the space is $\AS$.
As we saw in \autoref{ex:downhill}, the data is a function $\varphi \colon \RR^d \to \RR$ to minimize, as well as a set of $n$ starting point, and the algorithm $\Alg$ then produces a new set of $n$ points.

We say that two data points $d_1$ and $d_2$ in $\Cont(\AS) \times \AS^n$, respectively equal to $\varphi_1,X_1$ and $\varphi_2, X_2$ are \emph{related}, which we denote by
\begin{equation}
 \affrel{d_1}{d_2} 
\end{equation}
if
\begin{equation}
  \varphi_1 = \varphi_2 \circ \aff 
\qquad
 \text{and}
\qquad 
X_2 = \aff \mul X_1
,
\end{equation}
where the operation $\aff \mul X_1$ is the diagonal action.

The difference with \autoref{ex:downhill} is that the relation between $\varphi_1$ and $\varphi_2$ can no longer be expressed as $\varphi_2 = \varphi_1 \circ \aff\inv$, as $\aff$ is no longer required to be invertible.

The requirement that the algorithm $\Alg$ is affine equivariant in a stronger sense, is now that
\begin{equation}
  \affrel{d_1}{d_2} \implies \affrel{\Alg(d_1)}{\Alg(d_2)}
  .
\end{equation}
We will examine the consequences of such a stronger requirement in \autoref{ex:downhill3}, but we first put it in a formal setting.
\end{example}

\begin{remark}
  We use the word ``equivariance'', which is not completely correct.
  Indeed, equivariance is usually applied to natural transformations between functors from a group object.
  It may then perhaps be applied to any natural transformation.
  Howerver, as we shall see in \autoref{def:weaknattrans}, the algorithms which are fully affine equivariant fall short of being natural transformations, they are instead \emph{weak} natural tranformations (\autoref{def:weaknattrans}).
\end{remark}

We first describe the \emph{affine category}.
It consists of finite dimensional affine spaces as objects, and affine maps between affine spaces as morphisms.
To simplify the notations, we will identify all the affine spaces of the same dimensions, so the objects of the affine category are identified with the natural numbers:
\begin{equation}
\label{eq:identnataff}
  d \equiv \AS
  .
\end{equation}
As a result, $\Hom(m,n)$ denotes the affine maps between the affine spaces $\AS[m]$ and $\AS[n]$.

The data and computational domain are now replaced by the relevant functors.
These functors map an object in the affine cateogry (hence, a natural number) to an object in a category of smooth manifolds.

It turns out that we need the category of \emph{relations} associated to that of manifolds and smooth maps.
An object in that category is still a manifold, but a morphism between manifolds $\mathcal{M}$ and $\mathcal{N}$ is now a submanifold of $\mathcal{M} \times \mathcal{N}$.
If that submanifold is a graph, then this corresponds to a smooth map between $\mathcal{M}$ and $\mathcal{N}$, but this is otherwise a \emph{relation}.
We refer to~\cite{FrSc90} for a complete treatment of relation categories, also called \emph{allegories}.
For the general definition of full equivariance (\autoref{def:weaknattrans}), we only need to assume that $\Data $ and $\Comp$ are functors from the affine category to an allegory.

The data functor $\Data$ maps an affine space $d$ (identified with its dimension according to~\eqref{eq:identnataff}) to some manifold $\Data_d$.
A morphism $\aff \in \Hom(m,n)$ is mapped to a morphism $\Data(\aff) \in \Hom\paren[\big]{\Data_{m}, \Data_{n}}$ in the above allegory.

The data and computation objects are now indexed by an integer (the dimension).
We denote by $\affrel{x_1}{x_2}$ the fact that $x_1$ is related to $x_2$ by the affine map $\aff$, as in \autoref{ex:downhill2}.
This is also denoted by $(d_1,d_2) \in \Data(\aff)$.

The full equivariance condition is expressed as
\begin{equation}
\label{eq:fullaffequisimple}
\affrel{d_1}{d_2} \implies \affrel{\Alg_i(d_1)}{\Alg_j(d_2)}
,
\end{equation}
which can also be written as
\begin{equation}
(d_1,d_2)\in\Data(\aff) \implies \paren[\big]{\Alg(d_1),\Alg(d_2)} \in \Comp(\aff)
.
\end{equation}



What is the meaning of the full equivariance in the context of allegories?
The answer is that such a method is \emph{almost} a natural transformation between the functors $\Data$ and $\Comp$.
To understand this, we can reformulate condition~\eqref{eq:fullaffequisimple} using composition in the allegory.

But we must first address a small problem: the algorithm is a \emph{function}, and is thus not a morphism in the allegory (a relation).
But a function $\Alg$ naturally gives rise to a relation given by its graph, and we denote the corresponding relation by $\Rel{\Alg}$:
\begin{equation}
  \Rel{\Alg} \coloneqq \setc[\big]{(d,c)}{c = \Alg(d)}
  .
\end{equation}
We compute the composition
\begin{equation}
  \Rel{\Alg}_n \circ \Data(\aff) = \setc[\big]{(d_1,\Alg(d_2))}{(d_1,d_2) \in \Data(\aff)}
  ,
\end{equation}
and the composition
\begin{equation}
  \Comp(\aff) \circ \Rel{\Alg}_m = \setc[\big]{(d_1,c_2)}{(\Alg(d_1),c_2)\in\Comp(\aff)}
  .
\end{equation}
For an affine map $\aff \in \Hom(m,n)$, condition~\eqref{eq:fullaffequisimple} is thus
\begin{equation}
  \Rel{\Alg}_n \circ \Data(\aff) \subset \Comp(\aff) \circ \Rel{\Alg}_m
.
\end{equation}
Note that if the relation above was \emph{equality} instead of subset, it would be exactly the requirement that $F$ be a \emph{natural transformation} between the functors $\Data$ and $\Comp$.
This leads us to the definition of a weaker notion of a natural transformations in allegories.

\begin{definition}
\label{def:weaknattrans}
  Given two functors $\Data$ and $\Comp$ from a category $A$ to an allegory $B$, a \emph{weak natural transformation} is the data, for any object $M$ in the categry $A$, of a morphism $\Alg_M \in \Hom(\Data_M,\Comp_M)$, and such that for any object $M$ and $N$ in the category $A$, and any morphism $a\in\Hom(M,N)$ we have
\begin{equation}
\label{eq:weaknattrans}
  \Alg_N \circ \Data(a) \subset \Comp(a) \circ \Alg_M
.
\end{equation}
\end{definition}

The reader should compare~\eqref{eq:weaknattrans} with~\eqref{eq:nattrans}.


One can examine the meaning of full affine equivarience by breaking it into
particular cases.
Indeed, the affine category has two important subcategories: the category of \emph{injective} affine maps, and the category of \emph{surjective} affine maps.
We will call \emph{injective equivariance} and \emph{surjective equivariance} the property of being a weak natural transformation with respect to the corresponding subcategories.
For each dimension $d$, there is also a subcategory containing only the object $d$, and the invertible affine maps on that object: this is the category that we studied in \autoref{sec-equi}.
For each dimension, we will denote the corresponding equivariance by \emph{bijective equivariance}.

\begin{enumerate}
    \item \emph{Injective} equivariance generally means that if the data of the algorithm happens to lie in an affine subspace, then the result of the computation not only will lie on the subspace, but will also work exactly as if the lower dimensional version of the algorithms was used with that data.
    \item {Projective} equivariance generally indicates how the algorithm behaves with certain degenerate data. It highly depends on the algorithm. In the case of ODE integrators, it has a very understandable meaning (see \autoref{ex:bseries}).
    \item {Bijective} equivariance, is what we covered in \autoref{sec-equi}.
\end{enumerate}

\begin{example}[Polynomial interpolation and splines]
\label{ex:interpolation2}
  We revisit \autoref{ex:interpolation}.
  Now the dimension $d$ is arbitrary, and the algorithm works in any dimension.

The domain is $\AS^n$, the data of $n$ points in an affine space of dimension $d$.
The data functor $\Data$ maps the object $d$ of the affine category to $\Data_d = \AS^n$.
An affine map $\aff \in\Hom(m,n)$ is mapped to the relation $\Data(\aff)$ which we identify to the \emph{map} $\Data(\aff)(X) \coloneqq \aff \mul X$, where we used the diagonal action.
The computational functor $\Comp$ maps the object $d$ to $\Cont(\RR,\AS)$, and the relation $\Comp(\aff)$ is identified to the function $\Comp(\aff)(f) \coloneqq \aff \circ f$.

Note that in this case, both compositions occurring in \eqref{eq:weaknattrans} are \emph{graphs}, so the inclusion is in fact an equality, and $\Alg$ is in this case a \emph{natural transformation} between the functors $\Data$ and $\Comp$.

Injective equivariance here is related to a well known property of splines and interpolation: if the control points actually lie in a subspace, then the whole interpolating curve or spline, also lies in that subspace.
What is more, that curve is exactly the same as if the calculation had been done in a lower dimensional space instead.

Surjective equivariance is perhaps less intuitive: it means that interpolation commutes with affine projections.
Computing the interpolation of projected points on a smaller subspace gives the same result as projecting the interpolated curve instead.





\end{example}

\begin{example}[Downhill simplex]
\label{ex:downhill3}
  We now revisit \autoref{ex:downhill}.
  The data functor is $\Data_d = \Cont(\AS) \times \AS^n$ and a map $\aff\in\Hom(m,n)$ is mapped to $\Data(\aff)\in\Hom(\Data_m,\Data_n)$ defined as the relation
\begin{equation}
 \Data(\aff) = \setc{\paren[\big]{(\varphi_1,X_1),(\varphi_2,X_2)}}{\varphi_1 = \varphi_2 \circ \aff \qquad X_2 = \aff \mul X_1}
 .
\end{equation}
The computational functor is simply
\begin{equation}
  \Comp(\aff) = \setc{(X_1,X_2)}{X_2 = \aff \mul X_1}
  .
\end{equation}
Note that $\Comp(\aff)$ is in fact a graph, so it is associated to a function.

What does injective equivariance mean in that case?
We consider an \emph{injective} affine map $\aff \in \Hom(m,n)$.
The relation $\varphi_1 = \varphi_2 \circ \aff$ means that the function $\varphi_1 \in\Cont(\AS[m])$ to minimize is the \emph{restriction} of the function $\varphi_2(\AS[n])$, along the subspace given by the image of $\aff \in \Hom(m,n)$.
Equivariance means in this case is that: \emph{if one starts with a degenerate simplex}, i.e., if all the points lie in the subspace above, then the simplex algorithm will find the minimum \emph{in that subspace}, i.e., the minimum of the function $\varphi_1$.
\itodo{Picture here?}

Let us examine surjective equivariance.
We consider a \emph{surjective} affine map $\aff \in \Hom(m,n)$.
The relation $\varphi_1 = \varphi_2 \circ \aff$ now means that $\varphi_1$ is equal to $\varphi_2$ and is constant on the fibres (i.e., the level sets) of $\aff$.
This is an example of degenerate data.
What equivariance means in this case is that the simplex algorithms works \emph{fibrewise}, i.e., the result will not depend on where the initial points $X$ were chosen inside the fibres.
\end{example}

\begin{example}
\label{ex:bseries}
  We now look at the example that we understand perhaps best of all: numerical integration of ODEs.
  The data functor $\Data$ maps an affine space of dimension $d$ to $\Data_d = \mathcal{X}_0(\AS)$, the space of compactly supported vector fields on $\AS$.
  An affine map $\aff\in\Hom(m,n)$ is mapped to the relation $\Data(\aff)\in\Hom(\Data_m,\Data_n)$ defined by
  \begin{equation}
    \Data(\aff) = \setc{(f_1,f_2)}{f_2 \circ \aff = \Tan \aff \circ f_1}
    .
  \end{equation}
The computational functor $\Comp$ maps an object $d$ to $\Comp_d = \Diff(\AS)$. The relation $\Comp(\aff)\in\Hom(\Comp_m,\Comp_n)$ is defined by
\begin{equation}
  \Comp(\aff) = \setc{(\Phi_1,\Phi_2)}{\Phi_2 \circ \aff = \aff \circ \Phi_1}
  .
\end{equation}

The meaning of injective equivariance is known in numerical analysis as \emph{preservation of weak (affine) invariants}~\cite[\S\,IV.4]{HaLuWa06}.
An affine weak invariant is an affine subspace which is preserved by the flow of the vector field.
If that subspace is the image by an affine map $\aff$ of an affine space of smaller dimension (one can choose $\aff$ to be injective), the requirement of weak invariance is exactly that of being in relation with another vector field.
Injective equivariance thus means that: if a vector field has a weak invariant subspace, not only is it preserved by the numerical flow, but that numerical flow is the same as if computed in the lower dimensional subspace instead.

The meaning of surjective equivariance is particularly interesting.
Suppose that $\aff\in\Hom(m,n)$ is a surjective affine map.
The requirement that $(f_1,f_2) \in \Data(\aff)$ is that the flow of $f_1$ descends to the flow of $f_2$.
After change of variable, this can be rewritten as the differential equation
\begin{align}
  x' &= g_1(x)\\
  y' &= g_2(x,y)
  .
\end{align}
The property of surjective equivariance is that the numerical integrator behaves like the exact solution: the numerical flow descends to the  numerical flow of $f_2$.

We fully understand that case, as we can give a complete characterisation of the fully affine equivariant integrators in the sense above: these are exactly the integrators which have a \emph{B-Series}~\cite{bseries}.
\end{example}

\section*{Conclusion}

One of the biggest open questions is how weak natural transformations extend to other group actions.
Indeed, there are many examples of integrators on homogeneous spaces, which generalize the equivariance with respect to a \emph{group}.
However, as we saw, the equivariance with respect to the affine \emph{category} seems to be of the utmost importance.
We do not know of any other category for which a range of numerical algorithms are equivariant.

A particularly acute question is the characterization of numerical integrators on homogeneous spaces: as~\cite{aromatic} shows, group equivariance is not sufficient.
So what is the nonlinear equivalent of the affine category?
This is an area of ongoing research, but we speculate that they may be related to free Lie algebras, or possibly the related structures of post-Lie algebras.

\printbibliography%

\end{document}